\documentclass[11pt]{article}
\usepackage[a4paper]{anysize}\marginsize{3.5cm}{3.5cm}{1.3cm}{2cm}
\pdfpagewidth=\paperwidth \pdfpageheight=\paperheight
\usepackage{amsfonts,amssymb,amsthm,amsmath,eucal}
\usepackage{pgf}
\usepackage{bbm,array}
\usepackage{tikz} %Used for drawing picture in LaTeX
\usepackage{float}
\restylefloat{table}
\usepackage{subfigure}
\usepackage{caption}

\usetikzlibrary{arrows}
\theoremstyle{plain}
\newtheorem{thm}{Theorem}[section]
\newtheorem{theorem}[thm]{Theorem}

\newtheorem{lemma}[thm]{Lemma}
\newtheorem{proposition}[thm]{Proposition}

\theoremstyle{definition}

\newtheorem{remark}[thm]{Remark}
\newtheorem{example}[thm]{Example}

\newtheorem{problem}[thm]{Problem}

\newtheorem{thevarthm}[thm]{\varthmname}

\newenvironment{varthm*}[1]{\trivlist\item[]{\bf #1.}\it}{\endtrivlist}

\renewcommand\geq{\geqslant}

\newcommand\be{\begin{eqnarray*}}
\newcommand\ee{\end{eqnarray*}}

\newcommand\R{\mathbb R}

\newcommand\F{\mathbb F}

\renewcommand\P{\mathbb P}

\newcommand\calp{{\mathcal P}}
\newcommand\calr{{\mathcal R}}

\newcommand\call{{\mathcal L}}

\newcommand\calw{{\mathcal W}}

\newcommand\newop[2]{\def#1{\mathop{\rm #2}\nolimits}}
\newop\log{log}
\newop\ord{ord}
\newop\Gal{Gal}
\newop\SL{SL}
\newop\GL{GL}
\newop\Bl{Bl}
\newop\mult{mult}
\newop\mass{mass}
\newop\div{div}
\newop\codim{codim}
\newop\sing{sing}
\newop\vdim{vdim}
\newop\edim{edim}
\newop\Ass{Ass}
\newop\size{size}
\newop\reg{reg}
\newop\areg{areg}
\newop\asreg{asreg}
\newop\satdeg{satdeg}
\newop\supp{supp}
\newop\gin{gin}
\newop\ini{in}
\newop\vol{vol}
\newop\sat{sat}
\newop\length{length}
\newop\depth{depth}
\newop\characteristic{char}
\newcommand\eqnref[1]{(\ref{#1})}

\def\keywordname{{\bfseries Keywords}}%
\def\keywords#1{\par\addvspace\medskipamount{\rightskip=0pt plus1cm
\def\and{\ifhmode\unskip\nobreak\fi\ $\cdot$
}\noindent\keywordname\enspace\ignorespaces#1\par}}
\def\subclassname{{\bfseries Mathematics Subject Classification
(2000)}\enspace}
\def\subclass#1{\par\addvspace\medskipamount{\rightskip=0pt plus1cm
\def\and{\ifhmode\unskip\nobreak\fi\ $\cdot$
}\noindent\subclassname\ignorespaces#1\par}}

\newcommand\rounddown[1]{\left\lfloor#1\right\rfloor}

\definecolor{qqqqff}{rgb}{0,0,0}
\definecolor{uuuuuu}{rgb}{0,0,0}
\definecolor{zzttqq}{rgb}{0,0,0}
\definecolor{xdxdff}{rgb}{0,0,0}

\begin{document}

\author{A.~Czapli\'nski, M.~Dumnicki, \L .~Farnik, J.~Gwo\'zdziewicz, M.~Lampa-Baczy\'nska, G.~Malara,\\ T.~Szemberg\footnote{TS was partially supported by NCN grant UMO-2011/01/B/ST1/04875}, J.~Szpond, H.~Tutaj-Gasi\'nska}
\title{On the Sylvester--Gallai theorem for conics}
\date{\today}
\maketitle
\thispagestyle{empty}

\begin{abstract}
   In the present note we give a new proof of a result due to Wiseman and Wilson
   \cite{WisWil88} which establishes an analogue of the Sylvester-Gallai theorem
   valid for curves of degree two. The main ingredients of the proof come from
   algebraic geometry. Specifically, we use Cremona transformation of the projective
   plane and Hirzebruch inequality \eqnref{eq:Hirzebruch}.
\subclass{52C30, 05B30, 14Q10}
\end{abstract}

%*****************************************************************************

\section{Introduction}
\label{intro}
   Configurations of lines and points are a classical subject of study
   and a source of interesting results in combinatorics, geometry and algebra.
   One of the most celebrated results in this area is the Sylvester-Gallai Theorem.
\begin{theorem}[Sylvester-Gallai Theorem]\label{thm:SG}
   Let $\calp=\left\{P_1,\ldots,P_s\right\}$ be a finite number of points in the real projective plane.
   Then
   \begin{itemize}
   \item[a)] either all points are collinear;
   \item[b)] or there exists a line passing through exactly two points in the set $P_1,\ldots,P_s$.
   \end{itemize}
\end{theorem}
\begin{remark}
   The above result is of course also valid in the affine (euclidean) real plane.
   We have chosen the projective setting since it allows a particularly transparent
   proof of Theorem \ref{thm:Wiseman and Wilson}.
\end{remark}
   A line as in part b) of the above Theorem is called an \emph{ordinary line} for the set $\calp$.
   It is natural to wonder about the minimal number of ordinary lines in the dependence
   on the number of points $s$. Melchior \cite{Mel40} showed that there are at least $3$ such
   lines. It has been generalized to $\frac37s$ by Kelly and Moser, \cite[Theorem 3.6]{KelMos58}
   and further improved by Csima and Sawyer \cite[Theorem 2.15]{CsiSaw93}.
\begin{theorem}[Kelly and Moser, Csima and Sawyer]\label{thm:number of ordinary lines}
   For a set of $s$ non-collinear points in the real projective plane there are at least
   $$\frac37s\;\;\;\mbox{ ordinary lines.}$$
   Moreover, if $s\neq 7$, then the number of ordinary lines is at least $\frac{6}{13}s$.
\end{theorem}
   It has been conjectured by many authors that apart of the two cases constructed by
   Kelly and Moser, and Csima and Sawyer, the number of ordinary lines is bounded
   from below by $\frac12s$.
   Recently, in a ground breaking paper \cite{GreTao13}, Green and Tao proved that
   this is indeed the case for large values of $s$.

   There is a number of variations and generalizations of the Sylvester-Gallai Theorem,
   see e.g. \cite{BGS74}, \cite{BorMos90}, \cite{Gru09}.
   In most generalizations only linear objects are considered. This is in contrast
   with the following remarkable result proved in \cite{WisWil88} by Wiseman
   and Wilson.
\begin{theorem}[A Sylvester Theorem for conic sections]\label{thm:Wiseman and Wilson}
   Let $\mathcal{P}=\left\{P_1,\ldots,P_s\right\}$ be a finite number of points in the real projective plane.
   Then
   \begin{itemize}
   \item[a)] either all points lie on a conic;
   \item[b)] or there exists a conic $C$ passing through exactly five of the points in the set $\mathcal{P}$
      determined by these $5$ points (i.e. $C$ is the unique conic passing through these $5$ points).
   \end{itemize}
\end{theorem}
\begin{remark}
   A conic as in part b) of the above Theorem is called an \emph{ordinary conic} for the set $\calp$.
   It is irrelevant whether this conic is singular or not. In fact it might happen that \emph{all}
   ordinary conics for $\calp$ are singular, see Example \ref{ex:only singular}. The other extreme
   of all smooth ordinary conics is also possible, see Example \ref{ex:only smooth}.
\end{remark}
   The proof of this result presented in \cite{WisWil88} is quite involved. The purpose of the present note
   is twofold. First, we provide a simpler and more streamlined proof of Theorem \ref{thm:Wiseman and Wilson}.
   Second, it seems that the result of Wiseman and Wilson has not attracted as much attention as it deserves,
   we want to change this state of matters. In fact we find the result quite appealing and
   opening an unexplored path of research, with hight potential for substantial results.
   This fits well the philosophy presented in the recent survey by Tao \cite{Tao14}, to the effect that there
   are more hidden
   connections between various aspects of combinatorics and algebraic geometry.
   In section \ref{sec:further questions}
   we discuss some of natural further generalizations and pose some questions
   which hopefully will sparkle more interest and sparkle more research in this direction.
\section{Tools from algebraic geometry}
   The main tools we use in the proof of Theorem \ref{thm:Wiseman and Wilson} are the Cremona transformation
   and Hirzebruch inequality \eqnref{eq:Hirzebruch}.
   In this section we recall briefly these useful notions.

   We begin by the presentation of some basic properties of Cremona transformations.
   This part is valid over an arbitrary ground field.
   Let $F,G,H$ be non-collinear points in the projective plane $\P^2$. Let $h$ be the linear form
   defining the line determined by the points $F$ and $G$ and similarly: $g$ by $F,H$ and $f$ by $G$ and $H$.
   (By a slight abuse of notation we denote the lines by the same letters as their equations.)
   Then
   $$\P^2\ni(x:y:z) \stackrel{\varphi}
     {\longrightarrow} (g(x,y,z)\cdot h(x,y,z):f(x,y,z)\cdot h(x,y,z):f(x,y,z)\cdot g(x,y,z))\in\P^2$$
   is a birational automorphism of $\P^2$ (i.e. it is a $1:1$ map up to a codimension $1$ subvariety).
   It is the \emph{Cremona transformation
   based at the points} $F,G$ and $H$.
   After a projective change of coordinates, one may assume that the points
   $F$, $G$ and $H$ are the fundamental points (i.e. $(1:0:0)$, $(0:1:0)$ and $(0:0:1)$).
   Then the mapping $\varphi$ has a simple form
   $$\varphi:\P^2\ni(x:y:z)\mapsto (yz:xz:xy)\in\P^2.$$
   The planes before and after Cremona are schematically depicted in Figure \ref{fig: Cremona} below.
\begin{figure}[H]
\centering
   \begin{minipage}{0.4\textwidth}
   \centering
\begin{tikzpicture}[line cap=round,line join=round,x=1.0cm,y=1.0cm,scale=0.7]
\clip(0.64,-4.26) rectangle (8.78,3.42);
\draw [domain=0.86:8.] plot(\x,{(--10.5692-3.84*\x)/-1.58});
\draw [domain=0.86:8.] plot(\x,{(-7.8744-0.12*\x)/4.56});
\draw [domain=0.86:8.] plot(\x,{(-20.3948--3.96*\x)/-2.98});
\begin{scriptsize}
\draw [fill=qqqqff] (2.02,-1.78) circle (1.5pt);
\draw[color=qqqqff] (2.3,-1.48) node {$F$};
\draw[color=qqqqff] (5.15,0.5) node {$f$};
\draw[color=qqqqff] (2.6,0.5) node {$g$};
\draw[color=qqqqff] (4.5,-1.48) node {$h$};
\draw [fill=qqqqff] (6.58,-1.9) circle (1.5pt);
\draw[color=qqqqff] (6.82,-1.48) node {$G$};
\draw [fill=qqqqff] (3.6,2.06) circle (1.5pt);
\draw[color=qqqqff] (3.98,2.34) node {$H$};
\end{scriptsize}
\end{tikzpicture}
%   \caption{$ $ : $s=5$}\label{fig: s=5}
   \end{minipage}
   \quad
   \begin{minipage}{0.4\textwidth}
   \centering
\begin{tikzpicture}[line cap=round,line join=round,x=1.0cm,y=1.0cm,scale=0.7]
\clip(0.64,-4.26) rectangle (8.78,3.42);
\draw [domain=0.06:9.] plot(\x,{(--15.17-3.96*\x)/2.98});
\draw [domain=0.06:9.] plot(\x,{(-22.29--3.84*\x)/1.58});
\draw [domain=0.06:9.] plot(\x,{(-10.58--0.12*\x)/-4.56});
\begin{scriptsize}
\draw [fill=uuuuuu] (2.13,2.26474576271) circle (1.5pt);
\draw[color=uuuuuu] (2.38,2.68) node {$R_{FH}$};
\draw [fill=uuuuuu] (6.69,2.14) circle (1.5pt);
\draw[color=uuuuuu] (6.28,2.6) node {$R_{GH}$};
\draw [fill=uuuuuu] (5.11,-1.7) circle (1.5pt);
\draw[color=uuuuuu] (5.74,-1.5) node {$R_{FG}$};
\end{scriptsize}
\end{tikzpicture}
% \caption{$ $ : $s=6$}\label{fig: s=6}
\end{minipage}
\caption{$ $ : Cremona transformation}\label{fig: Cremona}
\end{figure}
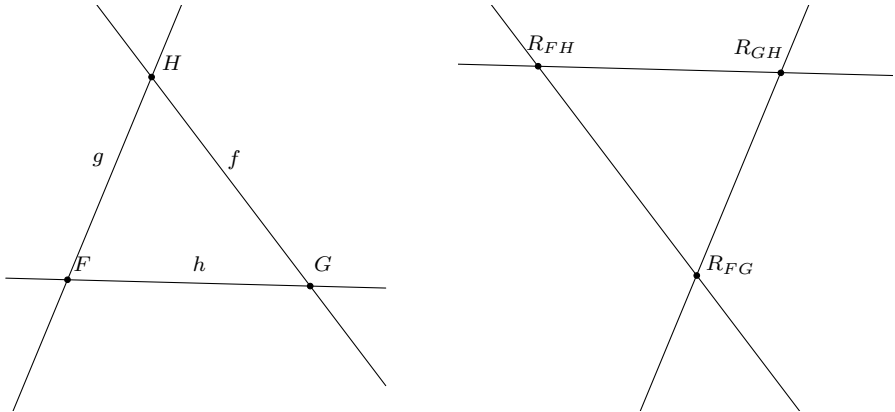
   The next Proposition collects basic properties of the Cremona map, which
   are relevant in the sequel. We refer to Dolgachev's masterpiece \cite{Dol12} for proofs
   and background. By a slight abuse of the notation the line defined by $f$ is denoted by $f$
   and similarly for $g$ and $h$.
\begin{proposition}\label{prop:Cremona}
   \begin{itemize}
      \item[a)] The Cremona map contracts the lines $f,g,h$ to points $R_{GH}$, $R_{FH}$, $R_{FG}$ respectively and it is $1:1$ away of them.
      \item[b)] The inverse mapping (in the category of birational maps) is also a Cremona transformation
         based at points $R_{FG}=\varphi(h)$, $R_{FH}=\varphi(g)$ and $R_{GH}=\varphi(f)$.
      \item[c)] Let $D$ be an irreducible curve, different from the three contracted lines,
         of degree $d$ with multiplicities $m_1$ at $F$, $m_2$ at $G$ and $m_3$ at $H$. Then its image $\varphi(D)$
         is an irreducible curve of degree $2d-m_1-m_2-m_3$ with multiplicities $d-m_2-m_3$ at
         $R_{GH}$, $d-m_1-m_3$ at $R_{FG}$ and $d-m_1-m_2$ at $R_{FG}$.
   \end{itemize}
\end{proposition}
   Now we pass to an inequality proved by Hirzebruch in the complex setting, see \cite[Section 3 and page 140]{Hir83}.
   The inequality itself is based on a very deep result due to Miyaoka, Yau and Bogomolov, \cite{Miy77}.
   Of course it remains valid for a configuration of real lines. Let $\call$ be an arrangement of $d$ lines.
   For $k\geq 2$, let $t_k(\call)$
   denote the number of points where exactly $k$ lines from
   $\call$ meet.
\begin{theorem}[Hirzebruch inequality]
   Let $\call$ be an arrangement of $d$ lines in the complex (or real) projective plane $\P^2$.
   Then
   \begin{equation}\label{eq:Hirzebruch}
      t_2(\call)+t_3(\call)\geq d+\sum\limits_{k\geq 5}(k-4)t_k(\call),
   \end{equation}
   provided $t_d=t_{d-1}=0$.
\end{theorem}
   In fact we will need the dual version of this inequality.
   To this end given a set of $s$ points in the projective plane let $t_i(\calp)$ denote
   the number of lines determined by this set (i.e. by pairs of points in the set),
   which pass through exactly $i$ points.
\begin{theorem}[Dual Hirzebruch inequality]
   Let $\calp$ be a set of $s$ distinct points. Assume that not all points are collinear and
   also not all but one point are collinear, then
   \begin{equation}\label{eq:dualHirzebruch}
      t_2(\calp)+t_3(\calp)\geq s+\sum\limits_{k\geq 5}(k-4)t_k(\calp),
   \end{equation}
\end{theorem}
   There is a similar inequality
   \begin{equation}\label{eq:Melchior}
      t_2(\calp)\geq 3+\sum\limits_{k\geq 4}(k-3)t_k(\calp)\;\;\mbox{ for }\; \calp\subset\P^2(\R),
   \end{equation}
   which was established by Melchior \cite{Mel40} using Euler
   formula applied to the partition of the \emph{real} projective plane
   given by the arrangement of lines. In the argument given
   in the next section one could work with this inequality
   instead of \eqnref{eq:dualHirzebruch}. However, in the view of Problem \ref{pro:complexSG}
   we prefer to work with a more general tool.
\section{Proof of Theorem \ref{thm:Wiseman and Wilson}}
   We begin by the following very useful observation.
\begin{lemma}[Main cases]\label{lem:Main cases}
   Let $\calp=\left\{P_1,\ldots,P_s\right\}$ be a finite set of points in the real projective plane $\P^2(\R)$.
   Then one of the following holds:
   \begin{itemize}
      \item[a)] all points in $\calp$ are collinear (i.e. $t_s(\calp)=1$) or;
      \item[b)] there is a line which contains exactly $3$ points from the set $\calp$ (i.e. $t_3(\calp)\geq 1$) or;
      \item[c)] there is a pair of ordinary lines intersecting in a point from $\calp$.
   \end{itemize}
\end{lemma}
\proof
   If a) holds, then we are done. So suppose that the points in $\calp$ are not collinear.
   If c) holds, then we are done again. So we are left with the situation when any two
   ordinary lines are disjoint (note that such configurations of points exist, see \cite[Figure 7]{BorMos90}).
   However then
   the number of ordinary lines is at most $\rounddown{\frac{s}{2}}$. Inequality \eqnref{eq:dualHirzebruch}
   implies then $t_3(\calp)\geq 1$.
\endproof
   Our proof of Theorem \ref{thm:Wiseman and Wilson} splits into three cases distinguished
   in Lemma \ref{lem:Main cases}. If the set $\calp$ consists of collinear points, then
   they are also contained in a conic and we are done. The next Lemma shows that the Theorem
   holds also in case b) of Lemma \ref{lem:Main cases}.
\begin{lemma}[Triple line]\label{lem:3 points}
   Let $\calp$ be a finite set in the real projective plane $\P^2(\R)$.
   If there is a line $L$ containing exactly $3$ points from $\calp$, then
   Theorem \ref{thm:Wiseman and Wilson} holds.
\end{lemma}
\proof
   Let $\calp'=\calp\setminus L$. If the set $\calp'$ is contained
   in a line $M$, then $\calp$ is contained in the union $L\cup M$,
   hence in a conic. Otherwise, there exists an ordinary line $M$ for $\calp'$.
   In that case, we take also $C=L\cup M$. There are exactly $5$ points
   from $\calp$ on $C$ and $C$ is uniquely determined by these points,
   since neither $L$, nor $M$ contains $4$ points from $\calp$. Note that
   it is irrelevant if the intersection point $L\cap M$ belongs to
   $\calp$ or not.
\endproof
   The rest of the proof deals with case c) of Lemma \ref{lem:Main cases}.
   The key idea here
   is to
   reduce the statement to Sylvester-Gallai theorem for lines
   applying the Cremona transformation based at the three points from $\calp$
   on the intersecting ordinary lines. The argument splits into several cases.

   Let $FG$ and $FH$ be ordinary lines for $\calp$ (intersecting in the point $F$).
   If their union contains the whole set $\calp$, then we are done. So we
   assume that this is not the case.
   We denote by $\calp''$ all points in $\calp$ contained in the union
   of the three lines determined by the points $F,G$ and $H$. In particular
   we have $F,G,H\in\calp''$. We call the
   residual set $\calp'$, i.e. $\calp'=\calp\setminus\calp''$.
   If $\calp'$ is empty, then $\calp$ is contained in the union of the line $GH$
   and any line through the point $F$, hence in a conic. So we assume that
   the set $\calp'$
   is non-empty.

   Let $\varphi$
   be the Cremona transformation based on the points $F,G$ and $H$
   and let $\calr=\varphi(\calp')$. In particular $\calr$ is non-empty.

   \textbf{Case 1.}
   We assume that all points in $\calr$ are collinear, contained in a line $L$. If the line
   $L$ is not uniquely determined, i.e. if there is just one point in $\calr$, then we take
   $L$ as in Subcase 1.b.

   \textbf{Subcase 1.a.}
   The line $L$ omits the points $R_{FG}, R_{FH}$ and $R_{GH}$ (so in particular there are
   at least two points in $\calr$). Then by Proposition \ref{prop:Cremona} c)
   the preimage of $L$ under $\varphi$ is a smooth conic $D$ passing through the points $F,G$ and $H$.
   The set $\calp$ is contained in the union of the four curves ($3$ lines and the conic $D$) indicated
   in the picture below.
\begin{figure}[H]
\centering
   \begin{minipage}{0.4\textwidth}
   \centering
\begin{tikzpicture}[line cap=round,line join=round,x=1.0cm,y=1.0cm,scale=0.7]
\clip(0.64,-4.26) rectangle (8.78,3.42);
\draw [domain=0.22:8.56] plot(\x,{(--10.5692-3.84*\x)/-1.58});
\draw [domain=0.22:8.56] plot(\x,{(-7.8744-0.12*\x)/4.56});
\draw [domain=0.22:8.56] plot(\x,{(-20.3948--3.96*\x)/-2.98});
\draw [rotate around={124.052497319:(4.12967764259,-0.680349313134)}] (4.12967764259,-0.680349313134) ellipse (2.89360787896cm and 2.37404756885cm);
\begin{scriptsize}
\draw [fill=qqqqff] (2.02,-1.78) circle (1.5pt);
\draw[color=qqqqff] (2.3,-1.48) node {$F$};
\draw [fill=qqqqff] (6.58,-1.9) circle (1.5pt);
\draw[color=qqqqff] (6.82,-1.56) node {$G$};
\draw [fill=qqqqff] (3.6,2.06) circle (1.5pt);
\draw[color=qqqqff] (3.88,2.34) node {$H$};
\draw [fill=qqqqff] (2.8,-2.72) circle (1.5pt);
\draw[color=qqqqff] (2.94,-2.44) node {$S$};
\draw [fill=qqqqff] (4.68,-3.42) circle (1.5pt);
\draw[color=qqqqff] (4.82,-3.14) node {$T$};
%\draw[color=black] (3.92,-2.9) node {$D$};
\draw [fill=xdxdff] (5.25716244605,-0.142135330999) circle (1.5pt);
\draw[color=xdxdff] (5.4,0.14) node {$U$};
%\draw [fill=qqqqff] (5.82,1.32) circle (1.5pt);
\draw[color=qqqqff] (5.96,1.6) node {$D$};
\end{scriptsize}
\end{tikzpicture}
%   \caption{$ $ : $s=5$}\label{fig: s=5}
   \end{minipage}
   \quad
   \begin{minipage}{0.4\textwidth}
   \centering
\begin{tikzpicture}[line cap=round,line join=round,x=1.0cm,y=1.0cm,scale=0.7]
\clip(0.64,-4.26) rectangle (8.78,3.42);
\draw [domain=-0.28:9.44] plot(\x,{(--15.1696474576-3.96*\x)/2.98});
\draw [domain=-0.28:9.44] plot(\x,{(-22.2872338983--3.84*\x)/1.58});
\draw [domain=-0.28:9.44] plot(\x,{(-10.5824135593--0.12*\x)/-4.56});
\draw [domain=2:9.44] plot(\x,{(-21.7892--4.42*\x)/6.6});
\begin{scriptsize}
\draw [fill=uuuuuu] (2.12644067797,2.26474576271) circle (1.5pt);
\draw[color=uuuuuu] (2.38,2.68) node {$R_{FH}$};
\draw [fill=uuuuuu] (6.68644067797,2.14474576271) circle (1.5pt);
\draw[color=uuuuuu] (6.28,2.6) node {$R_{GH}$};
\draw [fill=uuuuuu] (5.10644067797,-1.69525423729) circle (1.5pt);
\draw[color=uuuuuu] (5.74,-1.5) node {$R_{FG}$};
\draw [fill=xdxdff] (3.50350929689,-0.955104379965) circle (1.5pt);
\draw[color=xdxdff] (3.38,-0.5) node {$\varphi(S)$};
\draw [fill=xdxdff] (5.09349921707,0.109707051432) circle (1.5pt);
\draw[color=xdxdff] (5.02,0.5) node {$\varphi(T)$};
\draw[color=xdxdff] (7.82,1) node {$L=\varphi(D)$};
\end{scriptsize}
\end{tikzpicture}
% \caption{$ $ : $s=6$}\label{fig: s=6}
\end{minipage}
\caption{$ $ : Subcase 1.a}\label{fig: Subcase 1.a}
\end{figure}
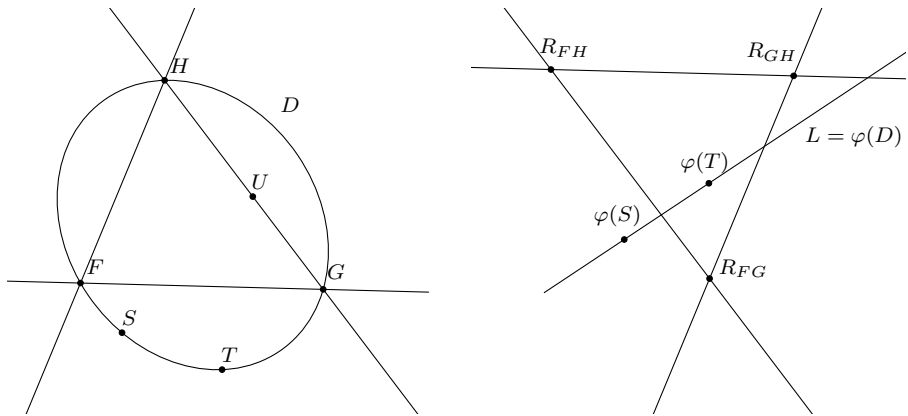
\definecolor{xdxdff}{rgb}{0,0,0}
\definecolor{qqqqff}{rgb}{0,0,0}
   If the line $GH$ is also an ordinary line for $\calp$, then
   all points in $\calp$ are on the conic $D$ and we are done.
   Otherwise, there
   is a point $U\in\calp$ on the line $GH$ as in the Figure \ref{fig: Subcase 1.a}.

   There are at least two points $S,T$ in $\calp'$. These points lie then on $D$.
   If the points $S,T,U$ are not collinear, then there is a single conic $C$
   determined by the points $S,T,U,F$ and $G$ and these are the only points in $\calp$ on $C$.
   So we are done.

   Note that $C$ is the union of lines if points $S,T,U$ are collinear and
   it is smooth otherwise.

%      In the remaining case (i.e. when $S,T,U$ are collinear) either there is another point $V\in\calp$ on the line $GH$
%      (which is then non-collinear with $S$ and $T$)
%      and we take the conic determined by $S,T,V,F$ and $G$, or there are just $3$
%      points from $\calp$ on the line $GH$ in which case we are done by Lemma \ref{lem:3 points}.

   \textbf{Subcase 1.b}
   The line $L$ goes through one of the points $R_{FG}, R_{FH}$ or $R_{GH}$.
   Note that $L$ cannot pass through any pair of these points, because
   the lines joining $R_{FG}, R_{FH}$ and $R_{GH}$ are not contained in the
   image of the Cremona transformation $\varphi$.

   We assume first that $L$ goes through $R_{GH}$.
   Then by Proposition \ref{prop:Cremona} c) the preimage $D$ of $L$ is a line
   passing through the point $F$ as indicated in the Figure below.
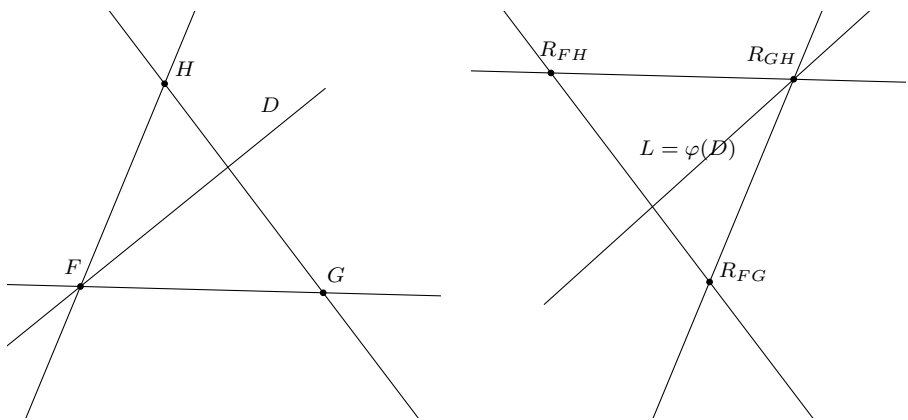
\begin{figure}[H]
\centering
   \begin{minipage}{0.4\textwidth}
   \centering
\begin{tikzpicture}[line cap=round,line join=round,x=1.0cm,y=1.0cm,scale=0.7]
\clip(0.64,-4.26) rectangle (8.78,3.42);
\draw [domain=0.42:9.62] plot(\x,{(--10.5692-3.84*\x)/-1.58});
\draw [domain=0.42:9.62] plot(\x,{(-7.8744-0.12*\x)/4.56});
\draw [domain=0.42:9.62] plot(\x,{(-20.3948--3.96*\x)/-2.98});
\draw [domain=0.42:6.62] plot(\x,{(-13.026--3.1*\x)/3.8});
\begin{scriptsize}
\draw [fill=qqqqff] (2.02,-1.78) circle (1.5pt);
\draw[color=qqqqff] (1.88,-1.4) node {$F$};
\draw [fill=qqqqff] (6.58,-1.9) circle (1.5pt);
\draw[color=qqqqff] (6.82,-1.56) node {$G$};
\draw [fill=qqqqff] (3.6,2.06) circle (1.5pt);
\draw[color=qqqqff] (3.98,2.34) node {$H$};
%\draw [fill=qqqqff] (5.82,1.32) circle (1.5pt);
\draw[color=qqqqff] (5.58,1.66) node {$D$};
\end{scriptsize}
\end{tikzpicture}
%   \caption{$ $ : $s=5$}\label{fig: s=5}
   \end{minipage}
   \quad
   \begin{minipage}{0.4\textwidth}
   \centering
\begin{tikzpicture}[line cap=round,line join=round,x=1.0cm,y=1.0cm,scale=0.7]
\clip(0.64,-4.26) rectangle (8.78,3.42);
\draw [domain=0.44:9.38] plot(\x,{(--15.1696474576-3.96*\x)/2.98});
\draw [domain=0.44:9.38] plot(\x,{(-22.2872338983--3.84*\x)/1.58});
\draw [domain=0.44:9.38] plot(\x,{(-10.5824135593--0.12*\x)/-4.56});
\draw [domain=2:9.38] plot(\x,{(-22.6936--5.24*\x)/5.76});
\begin{scriptsize}
\draw [fill=uuuuuu] (2.12644067797,2.26474576271) circle (1.5pt);
\draw[color=uuuuuu] (2.38,2.68) node {$R_{FH}$};
\draw [fill=uuuuuu] (6.68644067797,2.14474576271) circle (1.5pt);
\draw[color=uuuuuu] (6.28,2.6) node {$R_{GH}$};
\draw [fill=uuuuuu] (5.10644067797,-1.69525423729) circle (1.5pt);
\draw[color=uuuuuu] (5.74,-1.5) node {$R_{FG}$};
%\draw [fill=xdxdff] (4.7837811349,0.412050893553) circle (1.5pt);
\draw[color=xdxdff] (4.7,0.82) node {$L=\varphi(D)$};
\end{scriptsize}
\end{tikzpicture}
% \caption{$ $ : $s=6$}\label{fig: s=6}
\end{minipage}
\caption{$ $ : Subcase 1.b}\label{fig: Subcase 1.b}
\end{figure}
   In this situation $\calp$ is contained in the union of the line $D$ and the line $GH$.

   \textbf{Subcase 1.c}
   Now suppose that $L$ goes through the point $R_{FG}$ (the case $R_{FH}$ is analogous).
   Then its preimage is a line $D$ passing through $H$. If $GH$ is an ordinary line for
   $\calp$, then $\calp$ is contained in the union of $D$ and the line $FG$.
   If $GH$ contains exactly $3$ points from $\calp$, then we are done by Lemma \ref{lem:3 points}.
   In the remaining case $GH$ contains at least two points $U,V$ from $\calp$ distinct
   from the points $G$ and $H$. Also on $D$ there are at least two points $S,T$ from $\calp$
   distinct from the point $H$. Then the conic $C$ through $F, U, V, S$ and $T$ has these
   $5$ points in common with $\calp$ and it is determined by these points.
\begin{figure}[H]
\centering
   \begin{minipage}{0.4\textwidth}
   \centering
\begin{tikzpicture}[line cap=round,line join=round,x=1.0cm,y=1.0cm,scale=0.7]
\clip(0.64,-4.26) rectangle (8.78,3.42);
\draw [domain=0.14:9.48] plot(\x,{(--10.5692-3.84*\x)/-1.58});
\draw [domain=0.14:9.48] plot(\x,{(-7.8744-0.12*\x)/4.56});
\draw [domain=0.14:9.48] plot(\x,{(-20.3948--3.96*\x)/-2.98});
\draw [domain=0.14:9.48] plot(\x,{(--12.354-3.26*\x)/0.3});
\begin{scriptsize}
\draw [fill=qqqqff] (2.02,-1.78) circle (1.5pt);
\draw[color=qqqqff] (1.88,-1.4) node {$F$};
\draw [fill=qqqqff] (6.58,-1.9) circle (1.5pt);
\draw[color=qqqqff] (6.82,-1.56) node {$G$};
\draw [fill=qqqqff] (3.6,2.06) circle (1.5pt);
\draw[color=qqqqff] (3.98,2.34) node {$H$};
\draw [fill=qqqqff] (3.9,-1.2) circle (1.5pt);
\draw[color=qqqqff] (4.04,-0.92) node {$S$};
\draw [fill=xdxdff] (4.0675897589,-3.02114204673) circle (1.5pt);
\draw[color=xdxdff] (4.2,-2.74) node {$T$};
\draw [fill=xdxdff] (4.35522775018,1.0564087615) circle (1.5pt);
\draw[color=xdxdff] (4.5,1.34) node {$U$};
\draw [fill=xdxdff] (5.67806530413,-0.701455907499) circle (1.5pt);
\draw[color=xdxdff] (5.82,-0.42) node {$V$};
%\draw [fill=xdxdff] (3.78418302605,0.0585444502501) circle (1.5pt);
\draw[color=xdxdff] (4.04,0.26) node {$D$};
\end{scriptsize}
\end{tikzpicture}
%   \caption{$ $ : $s=5$}\label{fig: s=5}
   \end{minipage}
   \quad
   \begin{minipage}{0.4\textwidth}
   \centering
\begin{tikzpicture}[line cap=round,line join=round,x=1.0cm,y=1.0cm,scale=0.7]
\clip(0.64,-4.26) rectangle (8.78,3.42);
\draw [domain=0.24:9.] plot(\x,{(--15.1696474576-3.96*\x)/2.98});
\draw [domain=0.24:9.] plot(\x,{(-22.2872338983--3.84*\x)/1.58});
\draw [domain=0.24:9.] plot(\x,{(-10.5824135593--0.12*\x)/-4.56});
\draw [domain=0.24:9.] plot(\x,{(-10.2138386834--2.10730513084*\x)/-0.322659543068});
\begin{scriptsize}
\draw [fill=uuuuuu] (2.12644067797,2.26474576271) circle (1.5pt);
\draw[color=uuuuuu] (2.38,2.68) node {$R_{FH}$};
\draw [fill=uuuuuu] (6.68644067797,2.14474576271) circle (1.5pt);
\draw[color=uuuuuu] (6.28,2.6) node {$R_{GH}$};
\draw [fill=uuuuuu] (5.10644067797,-1.69525423729) circle (1.5pt);
\draw[color=uuuuuu] (5.74,-1.5) node {$R_{FG}$};
%\draw [fill=xdxdff] (4.7837811349,0.412050893553) circle (1.5pt);
\draw[color=xdxdff] (6.36,-3.16) node {$L=\varphi(D)$};
\draw [fill=xdxdff] (4.7,0.96) circle (1.5pt);
\draw[color=xdxdff] (5.34,1.24) node {$\varphi(S)$};
\draw [fill=xdxdff] (4.41773507035,2.80271550124) circle (1.5pt);
\draw[color=xdxdff] (5.06,3.08) node {$\varphi(T)$};
\end{scriptsize}
\end{tikzpicture}
% \caption{$ $ : $s=6$}\label{fig: s=6}
\end{minipage}
\caption{$ $ : Subcase 1.c}\label{fig: Subcase 1.c}
\end{figure}
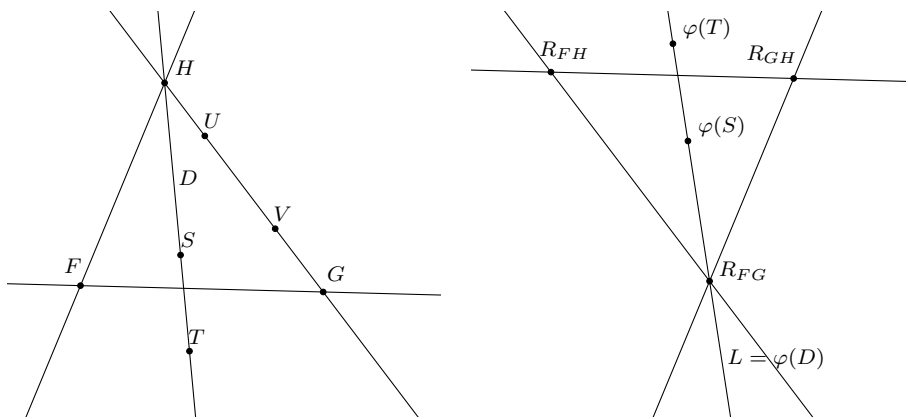

   \textbf{Case 2.}  We assume now that not all points in $\calr$ are collinear.
   Hence there exists an ordinary line $L$ for $\calr$. Let $\varphi(S)$ and $\varphi(T)$
   be the points in $\calr$ determining $L$.

   \textbf{Subcase 2.a.} We assume that $L$ does not pass through any of the
   points $R_{FG},R_{FH}$ and $R_{GH}$. This is the easiest case, since then
   Proposition \ref{prop:Cremona} c) implies that the preimage of $L$ under $\varphi$
   is an irreducible conic $C$ passing through points $S,T,F,G$ and $H$ and there are
   no more points from $\calp$ on $C$.

   \textbf{Subcase 2.b.} Now we assume that $L$ goes
   through the point $R_{FG}$ (the case when $L$ goes through $R_{FH}$
   is analogous). Then its preimage $D$ is a line passing
   through the point $H$. There are exactly $3$ points from $\calp$ on $D$,
   namely: $S,T$ and $H$. Note that the intersection point of $D$
   with the line $FG$ does not belong to $\calp$ (since $FG$ is an ordinary line).
   Thus we are done by Lemma \ref{lem:3 points}.
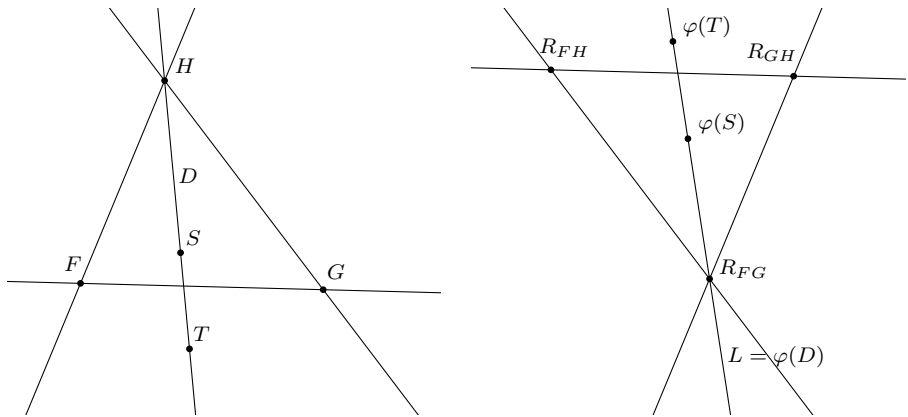
\begin{figure}[H]
\centering
   \begin{minipage}{0.4\textwidth}
   \centering
\begin{tikzpicture}[line cap=round,line join=round,x=1.0cm,y=1.0cm,scale=0.7]
\clip(0.64,-4.26) rectangle (8.78,3.42);
\draw [domain=0.14:9.48] plot(\x,{(--10.5692-3.84*\x)/-1.58});
\draw [domain=0.14:9.48] plot(\x,{(-7.8744-0.12*\x)/4.56});
\draw [domain=0.14:9.48] plot(\x,{(-20.3948--3.96*\x)/-2.98});
\draw [domain=0.14:9.48] plot(\x,{(--12.354-3.26*\x)/0.3});
\begin{scriptsize}
\draw [fill=qqqqff] (2.02,-1.78) circle (1.5pt);
\draw[color=qqqqff] (1.88,-1.4) node {$F$};
\draw [fill=qqqqff] (6.58,-1.9) circle (1.5pt);
\draw[color=qqqqff] (6.82,-1.56) node {$G$};
\draw [fill=qqqqff] (3.6,2.06) circle (1.5pt);
\draw[color=qqqqff] (3.98,2.34) node {$H$};
\draw [fill=qqqqff] (3.9,-1.2) circle (1.5pt);
\draw[color=qqqqff] (4.14,-0.92) node {$S$};
\draw [fill=xdxdff] (4.0675897589,-3.02114204673) circle (1.5pt);
\draw[color=xdxdff] (4.3,-2.74) node {$T$};
%\draw [fill=xdxdff] (3.78418302605,0.0585444502501) circle (1.5pt);
\draw[color=xdxdff] (4.04,0.26) node {$D$};
\end{scriptsize}
\end{tikzpicture}
%   \caption{$ $ : $s=5$}\label{fig: s=5}
   \end{minipage}
   \quad
   \begin{minipage}{0.4\textwidth}
   \centering
\begin{tikzpicture}[line cap=round,line join=round,x=1.0cm,y=1.0cm,scale=0.7]
\clip(0.64,-4.26) rectangle (8.78,3.42);
\draw [domain=0.24:9.] plot(\x,{(--15.1696474576-3.96*\x)/2.98});
\draw [domain=0.24:9.] plot(\x,{(-22.2872338983--3.84*\x)/1.58});
\draw [domain=0.24:9.] plot(\x,{(-10.5824135593--0.12*\x)/-4.56});
\draw [domain=0.24:9.] plot(\x,{(-10.2138386834--2.10730513084*\x)/-0.322659543068});
\begin{scriptsize}
\draw [fill=uuuuuu] (2.12644067797,2.26474576271) circle (1.5pt);
\draw[color=uuuuuu] (2.38,2.68) node {$R_{FH}$};
\draw [fill=uuuuuu] (6.68644067797,2.14474576271) circle (1.5pt);
\draw[color=uuuuuu] (6.28,2.6) node {$R_{GH}$};
\draw [fill=uuuuuu] (5.10644067797,-1.69525423729) circle (1.5pt);
\draw[color=uuuuuu] (5.74,-1.5) node {$R_{FG}$};
%\draw [fill=xdxdff] (4.7837811349,0.412050893553) circle (1.5pt);
\draw[color=xdxdff] (6.36,-3.16) node {$L=\varphi(D)$};
\draw [fill=xdxdff] (4.7,0.96) circle (1.5pt);
\draw[color=xdxdff] (5.34,1.24) node {$\varphi(S)$};
\draw [fill=xdxdff] (4.41773507035,2.80271550124) circle (1.5pt);
\draw[color=xdxdff] (5.06,3.08) node {$\varphi(T)$};
\end{scriptsize}
\end{tikzpicture}
% \caption{$ $ : $s=6$}\label{fig: s=6}
\end{minipage}
\caption{$ $ : Subcase 2.b}\label{fig: Subcase 2.b}
\end{figure}

   \textbf{Subcase 2.c.} Thus we are left with the situation that \emph{all}
   ordinary lines for $\calr$ go through the point $R_{GH}$. Note that
   this point is not contained in the set $\calr$. Let $t$ be the number
   of points in $\calr$. Then by Theorem \ref{thm:number of ordinary lines}
   there are at least $\frac37t$ ordinary lines for $\calr$. Each
   of these lines contains $2$ points from $\calr$, so that there
   are altogether at least $\frac67t$ points from $\calr$ on the union
   of these lines. Call this set $\calw$.
   Now, we consider the set $\calr'=\calr\cup\left\{R_{GH}\right\}$.
   All ordinary lines for $\calr'$ must be of the form: a line
   joining $R_{GH}$ with a point in $\calr\setminus\calw$. This
   implies that there are at most $\frac17t$ such lines, which
   contradicts Theorem \ref{thm:number of ordinary lines} for
   the set $\calr'$ consisting of $t+1$ points. Hence this subcase
   is not possible and the whole proof is finished.

\section{Examples and further questions}\label{sec:further questions}
   We begin with an example of a set of points $\calp$ such that
   \emph{every} ordinary conic for $\calp$ is singular. This shows
   that one cannot hope for Theorem \ref{thm:Wiseman and Wilson}
   to hold assuming $C$ smooth.
\begin{example}[Only singular ordinary conics]\label{ex:only singular}
   Let $C$ be a smooth conic and let $S$ be a point not on $C$. Let
   $L_1,L_2,L_3$ be three mutually distinct lines through $S$ intersecting
   $C$ in pairs of points $P_1,P_2$, $Q_1,Q_2$ and $R_1,R_2$ as indicated
   in the figure below.
\begin{figure}[H]
\centering
\begin{tikzpicture}[line cap=round,line join=round,>=triangle 45,x=1.0cm,y=1.0cm,scale=0.7]
\clip(-0.34,-3.96) rectangle (6.78,3.44);
\draw [rotate around={-0.32370152492:(3.07,-1.55)}] (3.07,-1.55) ellipse (2.51075041632cm and 1.78069302605cm);
\draw [domain=-0.34:6.78] plot(\x,{(--12.0386760811-5.67702001845*\x)/-1.55798660998});
\draw [domain=-0.34:6.78] plot(\x,{(--17.5126029995-6.15784380055*\x)/-0.121529585235});
\draw [domain=-0.34:6.78] plot(\x,{(--22.0298151576-5.39126155692*\x)/2.25181571919});
\begin{scriptsize}
\draw[color=black] (1.32,-0.48) node {$C$};
\draw [fill=qqqqff] (2.9,2.84) circle (1.5pt);
\draw[color=qqqqff] (3.34,3.12) node {$S$};
\draw [fill=xdxdff] (1.34201339002,-2.83702001845) circle (1.5pt);
\draw[color=xdxdff] (1.7,-2.52) node {$P_2$};
\draw[color=black] (4.,6.16) node {$a$};
\draw [fill=xdxdff] (2.77847041476,-3.31784380055) circle (1.5pt);
\draw[color=xdxdff] (3.16,-2.86) node {$Q_2$};
\draw[color=black] (3.18,6.16) node {$b$};
\draw [fill=xdxdff] (5.15181571919,-2.55126155692) circle (1.5pt);
\draw[color=xdxdff] (4.76,-2.2) node {$R_2$};
\draw[color=black] (1.74,6.16) node {$d$};
\draw [fill=uuuuuu] (2.1506830165,0.10962511214) circle (1.5pt);
\draw[color=uuuuuu] (1.92,0.54) node {$P_1$};
\draw [fill=uuuuuu] (2.84837878516,0.224378650118) circle (1.5pt);
\draw[color=uuuuuu] (3.12,0.64) node {$Q_1$};
\draw [fill=uuuuuu] (4.04999199216,0.0867073725263) circle (1.5pt);
\draw[color=uuuuuu] (4.18,0.46) node {$R_1$};
\end{scriptsize}
\end{tikzpicture}
\caption{$ $}\label{fig:only singular ordinary conics}
\end{figure}
   Then all ordinary conics for $\calp=\left\{S,P_1,P_2,Q_1,Q_2,R_1,R_2\right\}$
   split into a pair of lines
   through $S$. Indeed, an ordinary conic $C$ must pass through the point $S$
   and, since it contains altogether $5$ points from $\calp$, it must also pass
   through at least one pair
   of points on a line $L$ through $S$. But then $C$ and $L$ have at least
   $3$ points in common. By Bezout's Theorem, $L$ must be then a component of $C$,
   hence $C$ is a singular conic.
\end{example}
   On the other hand it might easily happen that \emph{all} ordinary conics for
   some set of points are \emph{smooth}.
\begin{example}[Only smooth ordinary conics]\label{ex:only smooth}
   Let $\calp$ be a set of points in \emph{general} position in the plane. Then
   there is an ordinary conic through any $5$ points from $\calp$ and
   all these conics are smooth (this is more or less the definition
   of the ''general position'').
\end{example}
   The Sylvester-Gallai Theorem fails in the finite characteristic. This is also
   the case for
   Theorem \ref{thm:Wiseman and Wilson}.
   The simplest counterexample is the following.
\begin{example}[Failure of the Theorem in finite characteristic]
   Let $\F$ denote the field with $5$ elements. Then $\P^2(\F)$
   consists of $31$ points. We consider the set $\calp$ consisting
   of all points in $\P^2(\F)$.
   Then any conic $C$ containing $5$ points from $\calp$ must contain at least
   one more point. Indeed, if $C$ is non-singular (and has $\geq 5$ points
   in $\P^2(\F)$), then it consists of exactly $6$ points. If it is singular,
   then it splits into two lines, each of them through $6$ points, so that
   there are altogether $11$ points from $\calp$ on $C$.
\end{example}
   The Sylvester-Gallai Theorem \ref{thm:SG} fails also over complex numbers.
   The simplest example is provided by the Hesse configuration, see \cite{ArtDol09}
   for details.

   We have expected that there exists also a complex counterexample
   to Theorem \ref{thm:Wiseman and Wilson}. However there are strong
   indications that this might not be the case. Of course, our proof
   of Theorem \ref{thm:Wiseman and Wilson} presented here, relies strongly
   on Theorem \ref{thm:SG}, so that it cannot be used in the complex case.
   It would be very interesting to know an answer to the following question.
\begin{problem}\label{pro:complexSG}
   Decide if Theorem \ref{thm:Wiseman and Wilson} is valid or not for points in the
   complex projective plane.
\end{problem}
   Once the problem is settled for curves of degree $1$ and $2$, it is natural
   to wonder what the situation is for curves of higher degree. Thus we repeat
   here the question which concludes article \cite{WisWil88}.
\begin{problem}[Curves of higher degree]\label{pro:higher degree}
   Let $\calp$ be a finite set of points in the projective plane
   and let $d$ be a positive integer.
   Does then one of the following hold:
   \begin{itemize}
   \item[a)] either $\calp$ is contained in a curve of degree $d$;
   \item[b)] or there exists a curve $C$ passing through exactly $\frac{(d+1)(d+2)}{2}-1$
   points in $\calp$ and determined by these points?
   \end{itemize}
\end{problem}
   This is not completely obvious if Problem \ref{pro:higher degree} is the right
   generalization of Theorems \ref{thm:SG} and \ref{thm:Wiseman and Wilson}. For
   example, for $d=3$ one might wonder instead if either $\calp$ is contained in a single
   cubic curve \emph{singular} in a point from $\calp$ or there exists such a curve
   determined by $\calp$.
\begin{remark}[Importance of the determined condition]
   Note that any line is determined by $2$ distinct points, so that it is not
   necessary to emphasize this condition in case b) of Theorem \ref{thm:SG}.
   This is no more the case for conics. In fact, it is very easy to show, that
   if not all points in a finite set $\calp$ are contained in a conic, then
   there exists a conic through exactly $5$ points in $\calp$.
   So that the critical point of Theorem \ref{thm:Wiseman and Wilson}
   is that there exist five points in $\calp$ which \emph{determine a single}
   conic.

\end{remark}
   Strangely enough the claim in the preceding Remark seems not easy to prove
   for curves of higher degree. So the following question can be viewed as the first
   step towards understanding Problem \ref{pro:higher degree}.
\begin{problem}
   Let $\calp$ be a finite set of points in the projective plane not contained in a curve
   of degree $d$. Is there a curve of degree $d$ passing through exactly $\frac{(d+1)(d+2)}{2}-1$
   points in $\calp$?
\end{problem}
\paragraph*{\emph{Acknowledgement.}}
   These notes originated during a workshop on Arrangements of Subvarieties held
   in Lanckorona in October 2014. We thank the Pedagogical University of Cracow
   for financial support.
%*****************************************************************************

%***************************************************************************** % Addresses

\bigskip \small

\bigskip
   Adam Czapli\'nski,
   Institut f\"{u}r Mathematik, Johannes Gutenberg-Universit\"at Mainz, D-55099 Mainz, Germany
   
\nopagebreak
   \textit{E-mail address:} \texttt{czaplins@uni-mainz.de}

\bigskip
   Marcin Dumnicki, \L ucja Farnik, Halszka Tutaj-Gasi\'nska,
   Jagiellonian University, Institute of Mathematics, {\L}ojasiewicza 6, PL-30-348 Krak\'ow, Poland

\nopagebreak
   \textit{E-mail address:} \texttt{Marcin.Dumnicki@im.uj.edu.pl}

   \textit{E-mail address:} \texttt{Lucja.Farnik@im.uj.edu.pl}

   \textit{E-mail address:} \texttt{Halszka.Tutaj@im.uj.edu.pl}

\bigskip
   Janusz~Gwo\'zdziewicz, Magdalena~Lampa-Baczy\'nska, Grzegorz~Malara, Tomasz Szemberg, Justyna Szpond,
   Instytut Matematyki UP,
   Podchor\c a\.zych 2,
   PL-30-084 Krak\'ow, Poland

\nopagebreak
%   \textit{E-mail address:} \texttt{@up.krakow.pl}

%   \textit{E-mail address:} \texttt{@up.krakow.pl}

  \textit{E-mail address:} \texttt{gwozd63@gmail.com}

  \textit{E-mail address:} \texttt{lampa.baczynska@wp.pl}

   \textit{E-mail address:} \texttt{gmalara@up.krakow.pl}

   \textit{E-mail address:} \texttt{szemberg@up.krakow.pl}

   \textit{E-mail address:} \texttt{szpond@up.krakow.pl}
%*****************************************************************************

\end{document}